\documentclass[11pt,reqno]{amsart}
\usepackage{amsfonts}
\usepackage{mathrsfs}
\usepackage{graphicx} 
\usepackage{epsfig}
\usepackage{amsmath, amsthm}
\usepackage{amssymb}
\usepackage[all]{xy}
\usepackage{color}

\numberwithin{equation}{section}

\theoremstyle{plain}

\numberwithin{equation}{section}

\theoremstyle{plain}
\newtheorem{thm}{Theorem}[section]
\newtheorem{prop}[thm]{Proposition}
\newtheorem{lem}[thm]{Lemma}
\newtheorem{cor}[thm]{Corollary}

\theoremstyle{definition}

\newtheorem{rem}[thm]{Remark}
\newtheorem{defn}[thm]{Definition}
\newtheorem{eg}[thm]{Example}
\newtheorem{subtitle}[thm]{}
\newtheorem{ex}{Exercise}[section]
\numberwithin{equation}{section}

\def\d{\delta}

\def\g{\gamma}
\def\G{\Gamma}
\def\K{\nabla}
\def\l{\lambda}

\def\n{\vert\/}
\def\o{\theta}

\def\ck{{\mathcal{K}}}
\def\cl{{\mathcal{L}}}
\def\cm{{\mathcal{M}}}

\def\cp{{\mathcal{P}}}

\def\li{\langle}
\def\ri{\rangle}
\def\n{|\/ }
\def\tr{{\rm tr}}
\def\bs{\bigskip}

\def\ss{\smallskip}

\def\ni{\noindent}
\def\ti{\tilde}
\def\p{\partial}

\def\Im{{\rm Im\/}}
\def\I{{\rm I\/}}

\def\diag{{\rm diag}}

\def\A{\mathbb{A}}
\def\C{\mathbb{C}}

\def\R{\mathbb{R} }

\newcommand{\beq}{\begin{equation}}
\newcommand{\eeq}{\end{equation}}
\newcommand{\beg}{\begin{eg}}
\newcommand{\eeg}{\end{eg}}
\newcommand{\bthm}{\begin{thm}}
\newcommand{\ethm}{\end{thm}}
\newcommand{\bprop}{\begin{prop}}
\newcommand{\eprop}{\end{prop}}
\newcommand{\bcor}{\begin{cor}}
\newcommand{\ecor}{\end{cor}}
\newcommand{\blem}{\begin{lem}}
\newcommand{\elem}{\end{lem}}
\newcommand{\bca}{\begin{cases}}
\newcommand{\eca}{\end{cases}}
\newcommand{\brem}{\begin{rem}}
\newcommand{\erem}{\end{rem}}
\newcommand{\bpm}{\begin{pmatrix}}
\newcommand{\epm}{\end{pmatrix}}
\newcommand{\bbm}{\begin{bmatrix}}
\newcommand{\ebm}{\end{bmatrix}}
\newcommand{\bvm}{\begin{vmatrix}}
\newcommand{\evm}{\end{vmatrix}}
\newcommand{\bdefn}{\begin{defn}}
\newcommand{\edefn}{\end{defn}}
\newcommand{\bsub}{\begin{subtitle}}
\newcommand{\esub}{\end{subtitle}}
\newcommand{\bex}{\begin{ex}}
\newcommand{\eex}{\end{ex}}
\newcommand{\ben}{\begin{enumerate}}
\newcommand{\een}{\end{enumerate}}

\def\det{{\rm det \/ }}

\def\rd{{\rm \/ d\/}}

\def\sech{{\rm sech\/}}

\def\an2{\hat A_{2n}^{(2)}}
\def\bn1{\hat B_n^{(1)}}

\def\B{\mathbb{B}}
\def\bbn1{\hat \B_n^{(1)}}
\def\A{\mathbb{A}}
\def\ban2{\hat \A_{2n}^{(2)}}

\def\fR{{\mathfrak R}}

\def\bh{\backslash}
\def\r0{\R^n \backslash \{0\}}
\def\bc{{\bf c}}

\def\bu{\bullet}
\def\n{\, | \,}

\begin{document}

\title{The Geometric Airy Curve Flow on $\R^n$} 

\author{Chuu-Lian Terng}
\address{Department of Mathematics\\
University of California at Irvine, Irvine, CA 92697-3875.  Email: cterng@math.uci.edu}

\maketitle

\centerline{Dedicated to Professor S. T. Yau on his 70th Birthday}

\section{Introduction}\label{aai}

The following third order curve flow in $\R^n$ was introduced in \cite{Ter15a},
\beq\label{ba}
\g_t= -\K_{e_1}^\perp H(\g(\cdot, t)),
\eeq
where $e_1= \g_x/||\g_x||$, $\K^\perp$ is the induced normal connection
and $H(\g(\cdot, t))$ is the mean curvature vector field of the curve $\g(\cdot, t)$.  
It is easy to see that
\beq\label{ba1} 
\g_t= -(\frac{1}{2} ||H||^2 e_1+ \K^\perp_{e_1} H(\g(\cdot, t)))
\eeq
is a curve flow in $\R^n$ that preserves arc-length parameter. Note that \eqref{ba1} and \eqref{ba} only differ by a tangent vector field, so they are essentially the same flows. In fact, the curve at time $t$ of \eqref{ba1} is the same curve  at time $t$ of \eqref{ba} reparametrized by the arc-length parameter.
Since the linear Airy equation is $$y_t= y_{xxx},$$
 we call \eqref{ba1} the {\it geometric Airy curve flow\/}. 
 
Curve flow \eqref{ba1} was considered by Langer and Perline in \cite{LP}. They proved the following:
\ben
\item[(i)] If $\g$ is a solution of \eqref{ba1} then there exists a parallel orthonormal frame $(e_2, \ldots, e_n)(\cdot, t)$ for the normal bundle of $\g(\cdot, t)$ with respect to $\K^\perp$ for each $\g(\cdot, t)$ such that the corresponding principal curvatures $k_i(\cdot, t)$ satisfy the {\it vector modified KdV equation\/} (vmKdV$_n$),
\beq\label{af2}
k_t= -(k_{xxx} + \frac{3}{2} ||k||^2 k_x),
\eeq
where $k=(k_1,\ldots, k_{n-1})^t$. 
\item[(ii)] There exist order $(2j-1)$ curve flow in $\R^n$ whose principal curvatures satisfy the $(2j-1)$-th vmKdV$_n$ flow. 
\een
 
Some main results of these paper are:
\ben 
\item We solve the Cauchy problem for the curve flow \eqref{ba1} with initial data $\g_0$ whose principal curvatures are rapidly decaying functions on $\R$. We also solve the periodic Cauchy problem for initial data with trivial normal holonomy. 
\item There is a Poisson structure on the space 
\beq\label{cj}
\cm_n=\{\g:S^1\to \R^n\n ||\g_x||\equiv 1\}
\eeq
so that the geometric Airy curve flow \eqref{ba1} is Hamiltonian.
\item  Equation \eqref{ba1} admits a sequence of commuting Hamiltonians, hence it has a sequence of contants of the motion. 
\item B\"acklund transformations (BT) are constructed, which generate a family of new solutions from a given solution of \eqref{ba1} by solving a linear Lax system.
\item An algorithm for constructing explicit soliton solutions of \eqref{ba1} is given.
\een

We mention some other integrable curve flows in flat spaces:  (i) Hasimoto showed that if $\g$ is a solution of the  {\it Vortex Filament equation\/} (VFE), $\g_t= \g_x\times \g_{xx}$,
 for curves in $\R^3$,
  then there exists a parallel frame $g(\cdot, t)$ along $\g(\cdot, t)$ so that $k_1+ik_2$ satisfies the NLS ( \cite{Has}). So we can apply soliton theory of the NLS to study geometric properties of the VFE (cf. \cite{LP91}, \cite{SY}, \cite{Ter15a}). 
 (ii) Pinkall considered in \cite{UP95}  the {\it central affine curve flow,\/} 
$\g_t= \g_{xxx} - 3q\g_x$,
 on $\R^2\bh 0$, where $x$ is the {\it central affine arc-length parameter\/} (i.e., $\det(\g,\g_x)=1$), and $q$ is the central affine curvature defined by $\g_{xx}= q \g$. He proved that if $\g$ is a solution of the central affine curve flow then its central affine curvature is a solution of the KdV, $q_t= q_{xxx} - 6 qq_x$. Its Hamiltonian aspect was given in \cite{CIM09}, and B\"acklund transformations were constructed in \cite{TWa}.
 The relation between central affine curve flows in $\R^n\bh 0$ and the Gelfand-Dickey flows on the space of n-th order linear differential operators on the line was given in \cite{CIM13} for $n=3$ and in \cite{TWb} for general $n$, and B\"acklund transformations were constructed in \cite{TerWu15a}. 
(iii) It was proved in \cite{TWc} that the {\it isotropic curve flow,\/}
$\g_t= \g_{xxx} - q \g_x$, in the light cone 
$\Sigma_{2,1}=\{x\in \R^{2,1}\n (x,x)=0\}$ preserves arc-length and its curvature $q$ is a solution of the KdV. For general $n$, let $\cm_{n+1,n}$ denote the space of curves $\g:\R\to \R^{n+1,n}$ such that the span of $\g, \ldots, \g_x^{(n-1)}$ is a Lagrangian subspace for all $x\in \R$ and $(\g_x^{(n)}, \g_x^{(n)})\equiv 1$. 
Isotropic curve flows were constructed in $\cm_{n+1,n}$ so that their isotropic curvatures were solutions of Drinfeld-Sokolov's KdV type flows associated to $B_n= o(n+1,n)$ and their Hamiltonian theory were given in  \cite{TWc}. B\"acklund transformations for these isotropic curve flows were constructed in \cite{TWd}.

It is interesting to note that the geometric Airy curve flow in $\R^n$, the isotropic curve flow in the null cone $\Sigma_{2,1}$ and the central affine curve flow  in $\R^2\bh 0$ are given by the re-parametrization of the third order flow $\g_t= \g_{xxx}$, where $x$ is the arc-length parameter for curves in $\R^n, \Sigma_{2,1}$ and is the central affine arc-length for curves in $\R^2\bh 0$.

This paper is organized as follows: We explain how to construct a parallel frame for a solution of \eqref{ba1} so that (i) hold and prove result (1) in section \ref{aag}. We review the Hamiltonian theory and construct B\"acklund transformations for the vmKdV$_n$ in section \ref{bz1}.    We prove (2) and (3) in section \ref{ay}. Results (4) to (5) are given in the last section.

\bs
\section{Geometric Airy Flow \eqref{ba1} on $\R^n$ and vmKdV$_n$}\label{aag}

In this section, we
\ben
\item[(i)] review the relation between the geometric Airy flow  \eqref{ba1} in $\R^n$ and the vmKdV$_n$ \eqref{af2}, 
\item[(ii)] solve the Cauchy problem on the line for \eqref{ba1} with initial data that has rapidly decaying principal curvatures,
\item[(iii)]  solve the periodic Cauchy problem for \eqref{ba1}. 
\een

Let $\g:\R\to \R^n$ be a smooth curve parametrized by its arc-length. 
 We call a smooth map $g=(e_1, e_2, \ldots, e_n):\R\to SO(n)$ a {\it parallel frame\/}  along $\g$ if $e_1=\g_x$ and $e_i$'s are parallel with respect to the induced normal connection for $2\leq i\leq n$. In other words, we have 
  \beq\label{ab1}
\bca (e_1)_x= \sum_{i=1}^{n-1} k_i e_{i+1},\\ (e_{i+1})_x= -k_i e_1, \quad 1\leq i\leq n-1.
\eca
\eeq
Equation \eqref{ab1} written in terms of $g$ is 
\beq\label{ab}
g^{-1}g_x= \bpm 0&-k^t\\ k&0\epm, \quad k=(k_1, \ldots, k_{n-1})^t.
\eeq
The function $k_i$ is the {\it principal curvature} of $\g$ with respect to the normal field $e_{i+1}$. 

When $n=2$,  a parallel frame for a plane curve is the Frenet frame and $k_1$ is the curvature of the plane curve. When $n>2$,  principal curvatures of a curve depend on the choice of parallel frames. In fact, if $g=(e_1, \ldots, e_n)$ is a parallel frames along $\g$ and $f=\diag(1, h)$ with a constant  $h\in SO(n-1)$ then $gf$ is also a parallel frame along $\g$ and the principal curvatures with respect to the new parallel frame $gf$ is $\ti k= h^{-1} {k}$, where $k= (k_1, \ldots, k_{n-1})^t$ and $\ti k
=(\ti k_1, \ldots, \ti k_{n-1})^t$. 

\blem\label{ada}
A curve flow in $\R^n$ of the form $\g_t= \sum_{i=1}^n A_i e_i$ preserves arc-length parameter if and only if 
\beq\label{adb}
e_1(A_1)= \sum_{i=2}^n A_i k_{i-1},
\eeq
where $(e_1(\cdot, t), \ldots, e_n(\cdot, t))$ is a parallel frame along $\g(\cdot, t)$ and $k_i(\cdot, t)$'s are the corresponding principal curvatures. 
\elem

\begin{proof}
Let $s(x,t)$ denote the arc-length parameter for $\g(\cdot, t)$, i.e., $s_x= ||\g_x||$. We use \eqref{ab1} to compute 
$$\frac{1}{2}\li \g_x, \g_x\ri_t= \li (\g_t)_x, \g_x\ri = ||\g_x||^2 \li (\g_t)_s, e_1\ri = ||\g_x||^2 ((A_1)_s -\sum_{i=2}^n A_i k_{i-1}).$$
This proves the Lemma. 
\end{proof}

The curve flow \eqref{ba} written in terms of parallel frame $(e_1,\ldots, e_n)(\cdot, t)$ is 
\beq\label{ac}
\g_t= -\sum_{i=1}^{n-1} (k_i)_s e_{i+1},
\eeq
where $k_i(\cdot, t)$ is the principal curvature of $\g(\cdot, t)$ with respect to the parallel normal field $e_{i+1}(\cdot, t)$ for $1\leq i\leq n-1$. It follows from Lemma \ref{ada} that 
\beq\label{ae3}
\g_t= -\frac{1}{2}\left( \sum_{i=1}^{n-1} k_i^2\right) e_1-\sum_{i=1}^{n-1} (k_i)_x e_{i+1},
\eeq
preserves the arc-length parameter.  Note that \eqref{ae3} is \eqref{ba1}. So we have the following:

\blem\label{ae1} The geometric Airy curve flow \eqref{ba1} in $\R^n$ preserves arc-length parameter and is \eqref{ae3} written in terms of parallel frame.
\elem

Recall an elementary lemma:

\blem\label{bb} Given smooth maps $A,B:\R^2\to so(n)$, then the following linear system 
\beq\label{bb1}
\bca g_x= gA,\\ g_t= gB, \eca
\eeq
is solvable for smooth $g:\R^2\to SO(n)$ if and only if 
\beq\label{bb2}
A_t= B_x+ [A, B].
\eeq
\elem

Next we give a proof of the following result in \cite{LP} and explain how to choose the suitable parallel frame (unique up to a constant in $SO(n-1)$) so that the principal curvature evolves according to vmKdV$_n$. 

\bthm\label{af} (\cite{LP})
Let $\g$ be a solution for \eqref{ba1}. Then there exists a  smooth map $g:\R^2\to SO(n)$ such that   $g(\cdot, t)$ is a parallel frame along $\g(\cdot, t)$ for each $t$ and $g$ satisfies
\beq\label{af1}
\bca g^{-1}g_x= \bpm 0& -k^t\\ k&0\epm,\\
g^{-1}g_t= \bpm 0 &-z^t\\ z& \xi\epm,
\eca
\eeq
where $k=(k_1, \ldots, k_{n-1})^t$,  and 
\begin{align*}
& z= -(k_{xx}+\frac{1}{2} ||k||^2 k), \\
&\xi=(\xi_{ij}), \quad \xi_{ij}= (k_i)_xk_j- k_i (k_j)_x, \quad 1\leq i, j\leq n-1.
\end{align*} 
Moreover, 
\ben
\item $k$ satisfies the vmKdV$_n$ \eqref{af2}.
\item If $\ti g:\R^2\to SO(n)$ satisfies \eqref{af1} with $\ti k$, then there exists a constant $c\in SO(n-1)$ such that $\ti g= g \diag(1, c)$ and $\ti k= c^{-1} k$. 
\een
\ethm

\begin{proof} Choose a smooth $h=(e_1, v_2, \ldots, v_n):\R^2\to SO(n)$  such that $h(\cdot, t)$ is a parallel frame along $\g(\cdot, t)$ for each $t$. Let $\mu_1, \ldots, \mu_{n-1}$ be the principal curvatures given by $h$. Then 
$$h^{-1}h_x= \bpm 0& -\mu^t\\ \mu & 0\epm, \quad \mu=(\mu_1, \ldots, \mu_{n-1})^t.$$
Write $h^{-1}h_t= \bpm 0& -z^t\\ z &C\epm$ with $z= (z_1, \ldots,  z_{n-1})^t$, and $C=(c_{ij})\in so(n-1)$. 
It follows from \eqref{ae3} and \eqref{ab} that we have
$$\li (e_1)_t, v_{i+1}\ri = \li (\g_t)_x, v_{i+1}\ri = -((\mu_i)_{xx} +\frac{1}{2}||\mu||^2 \mu_i).$$
Hence 
\beq\label{af3}
z_i= -((\mu_i)_{xx} +\frac{1}{2}||\mu||^2 \mu_i).
\eeq

Let $A:=h^{-1}h_x$ and $B=h^{-1}h_t$. By Lemma \ref{bb}, we have
$$A_t= B_x+[A,B].$$ Equate the $(22)$ block of the above equation to get 
$$(C_{ij})_x= -\mu_i(\mu_j)_{xx} + (\mu_i)_{xx} \mu_j = (-\mu_i (\mu_j)_x+ (\mu_i)_x\mu_j)_x.$$
Hence $C_{ij}(x,t)= (-\mu_i (\mu_j)_x+ (\mu_i)_x\mu_j)(x,t)+\eta(t)$ for some $\eta(t)\in o(n-1)$. 

Let $f:\R\to SO(n-1)$ be a solution of  $f_t= -\eta(t)f$, and 
$$g=(e_1, \ldots, e_n) := h\diag(1, f).$$ Then $g(\cdot, t)$ again a parallel frame along $\g(\cdot, t)$. The principal curvatures of $\g(\cdot, t)$ with respect to the new parallel frame $g$ is  $k(x,t)=f(t)^{-1} \mu(x,t)$.  A direct computation implies that $g$ satisfies \eqref{af1}. 

Let $P$ and $Q$ be the right hand sides of \eqref{af1}. By Lemma \ref{bb}, we have
\beq\label{af5} P_t= Q_x+ [P,Q].\eeq 
Equate the $(21)$ block of \eqref{af5} to see that $k$ satisfies \eqref{af2}. 

If $\ti g(\cdot, t)$ is another parallel frame along $\g(\cdot, t)$, then there exists $h:\R\to SO(n-1)$ such that $\ti g(x,t)= g(x,t) \diag(1, h(t))$ and the principal curvatures with respect to $\ti g$ is $\ti k(x,t)= h(t)^{-1} k(x,t)$. Note that 
\begin{align*}
\ti g^{-1}\ti g_t&= \diag(1,h^{-1}) g^{-1} g_t\diag(1, h)+ \bpm 0&0\\ 0& h^{-1}h_t\epm\\
& \quad=  \diag(1,h^{-1}) \bpm 0& -z(k)^t\\ z(k)& \xi(k)\epm \diag(1, h)+ \bpm 0&0\\ 0& h^{-1}h_t\epm.
\end{align*} But $\ti k(x,t)= h^{-1}(t) k(x,t)$ implies that 
$$z(\ti k)= h^{-1}z(k), \quad \xi(\ti k)= h^{-1}\xi(k) h.$$ Hence we have 
$$\ti g^{-1} \ti g_t= \bpm 0& -z(\ti k)^t\\ z(\ti k) & \xi(\ti k)\epm + \bpm 0&0\\ 0& h^{-1}h_t\epm.$$
By assumption $\ti g^{-1} \ti g_t=  \bpm 0& -z(\ti k)^t\\ z(\ti k) & \xi(\ti k)\epm$. So we have $h^{-1} h_t=0$, which implies that $h$ is a constant. This proves (2).
\end{proof}

The proof of Theorem \ref{af} also gives the following:

\bcor\label{ca} Given smooth $k:\R^2\to \R^{n-1}$,  then the linear system \eqref{af1} is solvable for $g:\R^2\to SO(n)$ if and only if $k$ satisfies \eqref{af2}. 
\ecor

Conversely, given a solution of vmKdV$_n$ we can construct a solution of the geometric Airy flow \eqref{ba}:

\bprop\label{cc}
Let $k$ be a solution of  vmKdV$_n$ \eqref{af2}. Then
\ben
\item there exists smooth $g=(e_1, \ldots, e_n):\R^2\to SO(n)$ satisfying \eqref{af1},
\item let 
\begin{align*}
Z&= -(\frac{||k||^2}{2} e_1 + \sum_{i=1}^{n-1} (k_i)_x e_{i+1}),\\
c(t)&= \int_0^t Z(0, t) \rd t.
\end{align*} Then 
\beq\label{ce}
\g(x,t)= c(t) + \int_0^x e_1(s,t) \rd s, 
\eeq 
is a solution of \eqref{ba1}.  
\een
\eprop

\begin{proof}
(1) follows from Corollary \ref{ca}.

The second equation of \eqref{af1} implies that 
$$(e_1)_t= -\sum_{i=1}^{n-1} ((k_i)_{xx} + \frac{1}{2} ||k||^2 k_i) e_{i+1}.$$
The first equation of \eqref{af2} is \eqref{ab1}. Use \eqref{ab1} and a direct computation to see that
$$ -\sum_{i=1}^{n-1} ((k_i)_{xx} + \frac{1}{2} ||k||^2 k_i) e_{i+1}= Z_x.$$
Since $c'(t)=Z(0, t)$, we have
$$\g_t=c'(t) +\int_0^x (e_1)_t \rd s=c'(t)+ \int_0^x Z_s\rd s= Z(x,t).$$
So $\g$ satisfies \eqref{ae3}, which is \eqref{ba1}.
\end{proof}

Let $\g:\R\to \R^n$ be a smooth curve parametrized by arc-length. If the principal curvatures with respect to a parallel frame $g$ are rapidly decaying  then the principal curvatures with respect to other parallel frame are also rapidly decaying. 

The Cauchy problem for the vmKdV$_n$ with rapidly decaying initial data $k^0:\R\to \R^n$ can be solved by the inverse scattering method (cf. \cite{BeaCoi85}). The periodic Cauchy problem for \eqref{af2} is solved in \cite{Bou}. Next we use these solutions to solve the Cauchy problems for the geometric Airy flow \eqref{ba1}. 

\bthm[{\bf Cauchy Problem for \eqref{ba1} on the line}]\label{ag} \

\ni Let $\g_0:\R\to \R^n$ be a smooth curve parametrized by its arc-length, $g_0(x)= (e_1^0(x), \ldots, e_n^0(x))$ a parallel frame along $\g_0$, and $k_1^0, \ldots, k_{n-1}^0$ the corresponding principal curvatures. Let $k:\R^2\to \R^{n-1}$ be the solution of the Cauchy problem 
\beq\label{ag2}
\bca k_t= -(k_{xxx} + \frac{3}{2} ||k||^2 k_x),\\ k(x,0)= k^0(x),\eca
\eeq
where $k^0=(k_1^0, \ldots, k^0_{n-1})$. 
Then there exists a unique $g=(e_1, \ldots, e_n):\R^2\to SO(n)$ satisfying \eqref{af1} with $g(0,0)= g^0(0)$.
Moreover, let $c(t)$ be as in Proposition \ref{cc}. Then
\beq\label{ag1}
\g(x,t)= \g_0(0) + c(t)+  \int_0^x e_1(s,t)\rd s
\eeq
is the solution of \eqref{ba1} with $\g(x,0)= \g_0(x)$. 
\ethm

\bthm[{\bf Periodic Cauchy Problem for \eqref{ba1}}]\label{ah} \

\ni Let $\g_0:\R\to \R^n$ be a smooth curve parametrized by its arc-length. Assume that $\g_0$ is periodic with period $2\pi$ and the holonomy of the induced normal connection along the closed curve $\g_0$ is trivial (so parallel frame $g^0$ and its principal curvature $k^0$ are periodic).  If $k(x,t)$ is the solution of \eqref{ag2} and is periodic in $x$, then we have the following:
\ben
\item Let $g=(e_1, \ldots, e_n):\R^2\to SO(n+1)$ be the solution of \eqref{af1} with $g(0,0)= g^0(0)$. Then  $g(x,t)$ is periodic in $x$ and $g(x,0)$ $=g^0(x)$. 
\item  The solution $\g$  of \eqref{ba1} defined by \eqref{ag1} is  periodic in $x$ and $\g(x,0)$ $= \g_0(x)$.  
\een
\ethm

\begin{proof} Corollary \ref{ca} gives the existence of $g$. Note that $g(x,0)$ and $g^0(x)$ satisfy the same linear ODE 
 and have the same initial data $y(0)= \g_0(0)$. Uniqueness of solutions of ODE implies that $g(x,0)= g^0(x)$.
 
 {\it Claim (A)}:  $g(x,t)$ is periodic in $x$. 
 
 To see this, we consider 
 $$y(t)= g(2\pi,t)-g(0,t).$$
 Recall that
\beq\label{ah1}
g_t= gB, \quad {\rm where\,\,} B= \bpm 0 & -z^t\\ z &\xi\epm,
\eeq
and $z=z(k)$ and $ \xi=\xi(k)$ are given as in Theorem \ref{af}. By assumption $k(x,t)$ is periodic in $x$. So we have 
\beq\label{ah2}
B(x+ 2\pi, t)= B(x,t)
\eeq 
for all $x,t\in \R$.
Use \eqref{ah1} and \eqref{ah2}  to see that
\begin{align*}
\frac{\rd y}{\rd t}&= g(2\pi, t) B(2\pi, t)- g(0,t) B(0,t)= g(2\pi,t)B(0,t)- g(0,t) B(0,t)\\
& \, = y(t) B(0,t). 
\end{align*}
Since $k$ is periodic in $x$,  $B(x,t)$ is periodic in $x$. Note that $g(x,0)= g^0(x)$ is periodic. So we have $y(0)=0$. Then the uniqueness of ODE implies that $y(t)\equiv 0$, i.e., $y$ is periodic. This proves Claim (A). 

{\it Claim (B)}: $\g(x,0)= \g_0(x)$.

We have proved $g(x,0)= g^0(x)$ is a parallel frame along $\g_0$. So $e_1(x,0)= (\g_0)_x$. By \eqref{ag1},  we have $\g_x= e_1$ and $\g(0,0)= \g_0(0)$. This shows that $\g(x, 0)=\g_0(x)$. 

{\it Claim (C)}:  $\g(x,t)$ defined by \eqref{ag1} is periodic in $x$. 

It follows from \eqref{ag1} that  
$$\eta(t)= \g(2\pi,t)-\g(0, t)= \oint e_1(s,t)\rd s.$$
Compute directly to see that
\begin{align*} 
\frac{\rd \eta}{\rd t} &= \oint (e_1)_t\rd s=\oint (\g_s)_t\rd s= \oint (\g_t)_s\rd s \\
&=-\oint \frac{\rd}{\rd s}\left(\frac{||k||^2}{2} e_1 + \sum_{i=1}^{n-1} (k_i)_s e_{i+1}\right) \rd s.
\end{align*}
Since $e_1(s,t), \ldots, e_n(s,t)$ and $k_i(s,t)$'s are periodic in $s$, we have  $\frac{\rd \eta}{\rd t}=0$. So $\eta(t)$ is a constant.   Note that $c(0)=0$,
$\eta(0)= \g(2\pi,0)-\g(0,0)= \g_0(2\pi)-\g_0(0)=0$. Hence $\eta(t)=0$ for all $t$ and Claim (C) is proved. 

Proposition \ref{cc} implies that $\g$ defined by \eqref{ag1} is the solution of \eqref{ba1} with $\g(x,0)= \g_0(x)$. 
\end{proof}

\bs
\section{The vmKdV hierarchy}\label{bz1}
 
 In this section, we review the construction of the vmKdV$_n$ hierarchy, explain its Hamiltonian theory, and construct  B\"acklund transformations. 

Let $\I_{n+1}$ denote the identity in $GL(n+1)$, and 
$$\I_{1,n}= \diag(1, -\I_n).$$
Let $\tau$ and $\sigma$ be involutions of $so(n+1,\C)$ defined by 
$$\tau(\xi)= \bar \xi, \quad \sigma(\xi)= \I_{1,n} \xi \I_{1,n}^{-1}.$$
Note that the fixed point set of $\tau$ is $so(n+1)$, $\tau\sigma= \sigma\tau$, and the $1, -1$-eigenspaces of $\sigma$ on $so(n+1)$  are
\begin{align*}
\ck&= \big\{\bpm 0&0\\ 0& y\epm\,\big|\, y\in so(n) \big\},\\
\cp&= \big\{\bpm 0& -z^t\\ z&0\epm \,\big|\,  z\in \R^{n\times 1}\big\}.
\end{align*}
Then we have $so(n+1)= \ck\oplus \cp$, and 
\beq\label{am}
[\ck,\ck]\subset \ck, \quad [\ck, \cp]\subset \cp, \quad [\cp, \cp]\subset \ck.
\eeq
The loop algebras we need to construct the vmKdV$_n$ hierarchy are: 
\begin{align*}
\cl&= \{\xi(\l)= \sum_{i\leq n_0} \xi_i \l^i\n n_0 \,\, {\rm some\,\, integer,\,\, } \xi_{2i}\in \ck, \xi_{2i-1}\in \cp\},\\
\cl_+&=\{\xi(\l)=\sum_{i\geq 0} \xi_i\l^i\in \cl\},\\
\cl_-&= \{\xi(\l)= \sum_{i<0} \xi_i\l^i \, \in \cl\}.
\end{align*} 

Note that 
$$\cl= \cl_+\oplus \cl_-$$
as linear subspaces. We call $(\cl_+, \cl_-)$ a {\it splitting\/} of $\cl$. 

\bdefn We say a map $\xi$ from $\C$ to $so(n+1,\C)$ or to $SO(n+1,\C)$ satisfies the {\it $(\tau, \sigma)$-reality condition\/} if
\beq\label{ar1}
\overline{\xi(\bar\l)} = \xi(\l), \quad \I_{1,n} \xi(-\l)\I_{1,n}^{-1}= \xi(\l).
\eeq
\edefn 

Note that $\xi(\l)=\sum_i\xi_i\l^i$ with $\xi_i\in so(n+1,\C)$ is in $\cl$ if and only $\xi$ satisfies the $(\tau,\sigma)$-reality condition \eqref{ar1}. 

A {\it vacuum sequence\/} is a linearly independent commuting sequence in $\cl_+$. Let 
\begin{align}
& a= e_{21}- e_{12}, \label{bp1}\\
& J_{2j-1}(\l)=  a\l^{2j-1}, \quad j\geq 1. \label{bp2}
\end{align}
Since $a\in \cp$, 
$\{J_{2j-1}\n j\geq 1\}$ is a vacuum sequence in $\cl_+$.

Next we use the standard method given in \cite{DS84}, \cite{TU00} to construct a soliton hierarchy from the splitting of $(\cl_+, \cl_-)$ of $\cl$ and the vacuum sequence $\{J_{2j-1}\n j\geq 1\}$.  

Given $z\in \R^{n-1}$, henceforth we will use the following notation: 
\beq\label{bs}
\Psi(z)= \bpm 0&0&0\\ 0&0& -z^t\\ 0& z&0\epm.
\eeq
Then
\beq\label{aj}
Y= [a\l, \cl_-]_+=[a, \cp]=\{\Psi(z)\n z\in \R^{(n-1)\times 1}\bigg\}.
\eeq

The following Theorem is known (cf. \cite{DS84}, \cite{TerUhl99b}). We include a proof here for completeness. 

\bthm \label{ak} Let $a$ be defined by \eqref{bp1}, and $Y$ by \eqref{aj}. Given $u\in C^\infty(\R, Y)$, then there exists a unique $Q(u,\l)= a\l + \sum_{i\geq 0} Q_i(u)\l^{-i}$
in $\cl$ satisfying 
\beq\label{ak1}
\bca [\p_x+ a\l + u, Q(u,\l)]=0, \\ Q(u,\l) \, {\rm is\, conjugate\, to \,} a\l.\eca
\eeq
Moreover, the $Q_i(u)$'s are differential polynomials of $u$ in $x$. 
\ethm

\begin{proof}
Since $Q_{2i}\in \ck$, $Q_{2i-1}\in \cp$, we can write 
 \begin{align}
 Q_{2j-1}&= \bpm 0& -y_{2j-1} & -\eta_{2j-1}^t\\ y_{2j-1}&0&0\\ \eta_{2j-1}&0&0\epm, \label{ake}\\
 Q_{2j}&= \bpm 0 & 0&0\\ 0&0& -z_{2j}^t\\ 0& z_{2j} & \xi_{2j}\epm, \label{ako}
 \end{align}
 where $y_{2j-1}:\R\to \R$, $\eta_{2j-1}, z_{2j}:\R\to \R^{n-1}$ and $\xi:\R\to so(n-1)$. 
 
Compare coefficients of $\l^{-i}$ of \eqref{ak1} to get the following recursive formula:
\beq\label{ak2}
(Q_i)_x+ [u, Q_i]+ [a, Q_{i+1}]=0.
\eeq
Assume that $u=\Psi(k)$  for some $k\in C^\infty(\R, \R^{n-1})$. Then $y_{2j-1}, \eta_{2j-1}, z_{2j}$ and $\xi_{2j}$ are differential polynomials of $k$.  We use \eqref{ak2} to see that 
\begin{align}
& (\eta_{2j-1})_x + y_{2j-1} k - z_{2j}=0,\label{ax1}\\
& (y_{2j-1})_x- k^t\eta_{2j-1}=0,\label{ax2}\\
& (z_{2j})_x - \xi_{2j} k+ \eta_{2j+1}=0, \label{ax3}\\
& (\xi_{2j})_x- kz_{2j}^t + z_{2j} k^t=0.\label{ax4}
\end{align}

We use the above equations, induction on $j$, and  direct computation to solve $Q_j(u)$ and see that entries of $Q_j(u)$ are differential polynomials of $u$.  Since $u=\Psi(k)$, $y_{2j-1}, \eta_{2j-1}, z_{2j}, \xi_{2j}$ are differential polynomials in $k$.
In particular, we obtain
\begin{align}
&z_0=k, \quad \xi_0=0,\label{ak3}\\
& y_1= -\frac{||k||^2}{2}, \quad \eta_1= -k_x, \label{ak4}\\
&  z_2= -(k_{xx} + \frac{1}{2}||k||^2 k),  \quad \xi_2= -kk_x^t +k_xk^t,  \label{ak5}\\
& y_3= k^tk_{xx} -\frac{||k_x||^2}{2} +\frac{3}{8} ||k||^4, \quad \eta_3= k_{xxx} +\frac{3}{2} ||k||^2 k_x, \label{ak6}
\end{align}
where $k=(k_1, \ldots, k_{n-1})^t$. 
\end{proof}

It follows from \eqref{ak2} that we have 
$$(Q_{2j-2}(u))_x+[u, Q_{2j-2}(u)]= [Q_{2j-1}(u), a]\in C^\infty(S^1, Y).$$ 
So the following is an evolution equation on $C^\infty(\R, Y)$,
\beq\label{al}
u_t= [\p_x+ u, Q_{2j-2}(u)] = [Q_{2j-1}(u), a].
\eeq
Equation \eqref{al} for $u$ written in terms of $k$ is
\beq\label{al1}
k_t= (z_{2j-2}(k))_x-\xi_{2j-2}(k)k,
\eeq
where $u=\Psi(k)$, and $z_{2j-2}(k)$ and $\xi_{2j-2}(k)$ are the $(32)$ and $(33)$ blocks of $Q_{2j-2}(u)$ respectively.  
Note that for $j=2$,  \eqref{al1} is $k_t= (z_2)_x-\xi_2 k$, where $z_2, \xi_2$ are defined by \eqref{ak5}. So it is the vmKdV$_n$ \eqref{af2}.

\bdefn
We call \eqref{al} (or \eqref{al1}) the {\it $(2j-1)$-th vmKdV$_n$ flow} and this sequence of flows the {\it vmKdV$_n$ hierarchy\/}. 
\edefn

It follows from \eqref{ak2} that the coefficients of $\l^i$ with $i>0$ of 
$$[\p_x+ a\l+ u, (Q(u,\l)\l^{2j-2})_+]$$
 are zero. So we have the following well-known existence of Lax pair.  

\bprop \label{cg} The following statements for smooth $u:\R^2\to Y$ are equivalent:
\ben
\item $u$ is a solution of the $(2j-1)$-th vmKdV$_n$ flow \eqref{al}.
\item The following linear system is solvable for $g:\R^2\to SO(n)$
\beq\label{cf}
\bca g_x= g u,\\ g_t= gQ_{2j-2}(u).\eca
\eeq
\item $u$ satisfies 
\beq\label{an}
u_t= [\p_x+ a\l+ u, (Q(u,\l)\l^{2j-2})_+],
\eeq
\item The following linear system is solvable for $E(x,t,\l)\in O(n+1,\C)$
\beq\label{an1}
\bca E^{-1}E_x= a\l+ u, \\ E^{-1}E_t= (Q(u,\l)\l^{2j-2})_+,\\
\overline{E(x,t,\bar\l)}= E(x,t,\l), \quad \I_{1,n}E(x,t,-\l) \I_{1,n}^{-1} = E(x,t,\l).
\eca
\eeq
\een
\eprop
 
 \bdefn We call \eqref{an1} a {\it Lax system\/} for \eqref{al}, and a solution $E(x,t,\l)$ of \eqref{an1}  a {\it frame\/} of a solution $u$ of \eqref{al} provided $E(x,t,\l)$ is holomorphic for all $\l\in \C$.  
 \edefn
 
 \ss
 \bsub {\bf B\"acklund transformations for vmKdV$_n$}\
 
We use the loop group factorization method given in \cite{TU00} to construct B\"acklund transformations for the vmKdV$_n$ hierarchy. 
 Let $L_+$ denote the group of holomorphic maps $f:\C\to SO(n+1, \C)$ that satisfies the $(\tau,\sigma)$-reality condition \eqref{ar1}, and $\fR_-$ the group of rational maps $f:\C\cup \{\infty\}\to SO(n+1,\C)$ that satisfies the $(\tau, \sigma)$-reality condition \eqref{ar1} and $f(\infty)=\I_{n+1}$.  
 
 First we recall the following results in \cite{TU00}. 
 
 \bthm\label{at} (\cite{TU00})
\ben
\item Given $g\in \fR_-$ and $f\in L_+$, then there exists unique $\ti g\in \fR_-$ and $\ti f\in L_+$ such that 
$gf= \ti f \ti g$.
\item Let $E(x,t,\l)$ be a frame of a solution $u$ of the $(2j-1)$-th vmKdV$_n$ flow \eqref{al}, $g\in \fR_-$,  $\ti E(x,t,\cdot)\in L_+$, and $\ti g(x,t,\cdot)\in \fR_+$ satisfying
$$g(\l) E(x,t,\l) = \ti E(x,t,\l) \ti g(x,t,\l).$$
Expand 
$$\ti g(x,t,\l)= \I+ \ti g_{-1}(x,t)\l^{-1} + \ldots.$$
Then 
\beq\label{at1}
\ti u= u+[a, \ti g_{-1}]
\eeq
is again a solution of \eqref{al} and $\ti E(x,t,\l)$ is a frame of $\ti u$. 
\een
  \ethm

We call an element $f$ in $\fR_-$ {\it a simple element\/} if $f$ cannot be written as a product of two elements in $f\R_-$. Given a simple element $g\in \fR_-$ and an $f\in L_+$, if we can  write down an explicit formula for the factorization of $gf= \ti f\ti g$ with $\ti f\in \fR_-$ and $\ti g\in L_+$, then $u\mapsto \ti u$ of Theorem \ref{at} (2) gives a B\"acklund transformation for \eqref{al}. 

 Let $s\in \R\bh 0$, $\pi$ a Hermitian projection of $\C^{n+1}$ satisfying 
\beq\label{au}
\bar \pi= \I_{1,n} \pi \I_{1,n}^{-1}, \quad \pi \bar \pi=\bar\pi \pi=0,
\eeq
or equivalently,
\beq\label{au1}
\bar V= JV, \quad \li V, JV\ri=0.
\eeq

\brem\label{az3} $v=\bpm y_0\\ y_1\epm \in \C^{n+1}$ with $y_0\in \C$ and $y_1\in \C^n$ satisfying 
\beq\label{az}
\bar v= \I_{1,n} v, \quad v^t v=0
\eeq
if and only if $y_0\in \R$, $y_1= i\bc$ for some $\bc\in \R^n$ and $||\bc||=|y_0|$. In particular, if $v\not=0$ then $y_0\not=0$.  Note that $v$ satisfies \eqref{az} if and only if the Hermitian projection $\pi$ onto $\C v$ satisfies \eqref{au}. 
\erem

\blem\label{az1}
Given $f\in L_+$ and a non-zero $v\in \C^{n+1}$ satisfying \eqref{az}, let $\ti v:= f(-is)^{-1}v$. Then $\ti v$ satisfies \eqref{az}, $\ti v= \bpm \ti y_0\\ i\ti \bc\epm$ for some $y_0\in \R\bh 0$, $\bc\in \R^n$, and $||\bc||= |\ti y_0|$. 
\elem

\begin{proof} Note that if $\in L_+$ then $f:\C\to SO(n+1,\C)$. So we have
\beq\label{ar2}
f(\bar\l)^*f(\l)=\I_{n+1}.
\eeq
  Use the reality condition \eqref{ar1}, \eqref{ar2} and \eqref{az} to compute 
$$\I_{1,n} \overline{\ti v}= \I_{1,n}\overline{f(-is)^{-1}}\bar v= \I_{1,n} f(is)^t \bar v= f(-is)^t \I_{1,n}\bar v= f(-is)^{-1} v= \ti v.$$
Similarly,
\begin{align*}
\li \overline{\ti v}, \ti v\ri &= \li \I_{1,n}\ti v, \ti v\ri = \li \I_{1,n}f(-is)^{-1} v, \ti v\ri= \li \I_{1,n} f(is)^\ast v, \ti v\ri \\ &= \li v, f(is)\I_{1,n} \ti v\ri = \li v, \I_{1,n}f(-is)\ti v\ri = \li v, \I_{1,n} v\ri= \li v, \bar v\ri=0.
\end{align*}
The rest of the Lemma follows from Remark \ref{az3}.
\end{proof}

The following Proposition gives simple elements in $\fR_-$ and a permutability formula (i.e., a relation among simple elements) for $\fR_-$. This type result was first proved for the the group of rational maps $f:\C\to GL(n,\C)$ that satisfies \eqref{ar2}  in \cite{TU00}. 

\bprop\label{av1} (\cite{DFG})
Let $s\in \R\bh 0$ and $\pi$ a Hermitian projection of $\C^{n+1}$ satisfying \eqref{au}. Then 
\beq\label{av}
\phi_{is, \pi}= (\I+ \frac{2is}{\l-is}\bar\pi^\perp) (\I-\frac{2is}{\l+is} \pi^\perp)= \I+ \frac{2is}{\l-is} \pi - \frac{2is}{\l+is} \bar \pi.
\eeq is in $\fR_-$. Moreover, we have the following:
\ben
\item  $\phi_{is,\pi}^{-1}=\phi_{-is, \pi}$.
\item Given $f\in L_+$, let let $\ti V= f(-is)^{-1}(\Im\pi)$, and $\ti\pi$ the Hermitian projection onto $\ti V$. Then $\ti\pi$ satisfies \eqref{au} and
\beq\label{av2}
\ti f:= \phi_{is,\pi}f \phi^{-1}_{is, \ti\pi}
\eeq
is in $L_+$. In other words, we have factored
$\phi_{is,\pi} f= \ti f \phi_{is, \ti\pi}$ with $\ti f\in L_+$ and $\phi_{is, \ti\pi}\in \fR_-$. 
\item Let $s_1^2\not= s_2^2\in \R\bh 0$, $\pi_1, \pi_2$ Hermitian projections of $\C^{n+1}$ satisfying \eqref{au}, and $\tau_1, \tau_2$ Hermitian projections onto $\C^{n+1}$. Then 
\beq\label{av4}
\phi_{is_2,\tau_2}\phi_{is_1, \pi_1} = \phi_{is_1, \tau_1} \phi_{is_2, \pi_2}
\eeq 
if and only if
\beq\label{av3}
\bca 
\Im\tau_1= \phi_{-is_2, \pi_2}(-is_1)(\Im \pi_1),\\
\Im\tau_2= \phi_{-is_1, \pi_1}(-is_2)(\Im \pi_2).
\eca
\eeq
\een
\eprop

As a consequence of Theorem \ref{at} and Proposition \ref{av1}(2), we obtain BTs for the $(2j-1)$-th flow \eqref{al}.

\bthm[{\bf B\"cklund Transformation for \eqref{al}}]\label{aw} \

\ni Let $s\in \R\bh 0$,  $\pi$ a Hermitian projection of $\C^{n+1}$ satisfying \eqref{au}, and $\phi_{is,\pi}$ defined by \eqref{av}.  Let $u$ be a solution of the $(2j-1)$-th vmKdV$_n$ flow \eqref{al}, and $E(x,t,\l)$ a  frame of $u$. Let $\ti\pi(x,t)$ be the Hermitian projection of $\C^{n+1}$ onto 
$$\ti V(x,t)= E(x,t,-is)^{-1}(\Im \pi).$$
Then 
\beq\label{aw1}
\ti u= u+ 2is [a, \ti\pi-\overline{\ti \pi}]
\eeq is a new solution of \eqref{al} and
\beq\label{aw2}
\ti E(x,t,\l)= \phi_{is,\pi}(\l) E(x,t,\l) \phi_{is, \ti\pi(x,t)}^{-1}
\eeq
is a frame of $\ti u$. 
\ethm

\esub

\bsub {\bf Permutability for BTs of \eqref{al}}\

First we recall the following result. 

  \bthm\label{bc} (\cite{TU00}) Let $E(x,t,\l)$ be the frame of a solution $u$ of the $(2j-1)$-th vmKdV$_n$ flow \eqref{al} such that $E(0,0,\I)=\I_{n+1}$.  Let $f\in \fR_-$, and $\ti u$ the solution constructed in Theorem \ref{aw} from $E$ and $f$. Then $f \bullet u= \ti u$ defines an action of $\fR_-$ on the space of solutions of \eqref{al}. 
\ethm

\bthm\label{bq}
 Let $s_i, \pi_i, \tau_i$ with $i=1,2$ be as in Proposition \ref{av1}(3), and $E$ the frame of a solution $u$ of \eqref{al} with $E(0,0,\l)=\I_{n+1}$. Let $\ti\pi_1, \ti \pi_2, \ti\tau_1, \ti\tau_2$ be the Hermitian projection of $\C^{n+1}$ such that 
\begin{align*}
\Im (\ti\pi_i(x,t))&= E(x,t,-is_i)^{-1}(\Im\pi_i), \quad i=1,2,\\
\Im(\ti\tau_1(x,t) )&= \phi_{-is_2, \ti\pi_2(x,t)}(\Im\ti\pi_1(x,t)),\\
\Im(\ti\tau_2(x,t))&= \phi_{-is_1, \ti\pi_1(x,t)}(\Im\ti\pi_2(x,t)).
\end{align*}
Then 
$$u_{12} = u_1+ 2is_2[a, \ti \tau_2-\overline{\ti\tau_2}], \quad u_{21} =  u_2+ 2is_1 [a, \ti \tau_1-\overline{\ti\tau_1}]$$
are solutions of \eqref{al} and $u_{12}= u_{21}$, where $u_i= u+ 2is_i[a, \ti\pi_i-\overline{\ti \pi_i}]$ for $i=1,2$. Moreover, 
$$E_{12}= \phi_{is_2, \tau_2}\phi_{is_1, \pi_1} E \phi_{is_1, \ti\pi_1}^{-1} \phi_{is_2, \ti\tau_2}$$ is the frame of $u_{12}$ with  $E_{12}(0,0,\l)= \I_{n+1}$. 
\ethm 

\begin{proof} Let $u_j= \phi_{is_j,\pi_j}\bu u$. Since $\bu$ is an action, we have
\begin{align*}
u_{12}&== \phi_{is_2,\tau_2}\bu u_1= \phi_{is_2,\tau_2}\bu (\phi_{is_1,\pi_1}\bu u)= (\phi_{is_2,\tau_2}\phi_{is_1,\pi_1})\bu u,\\
u_{21}&= \phi_{is_1,\tau_1}\bu (\phi_{is_1,\pi_2}\bu u)= \phi_{is_1,\tau_1}\bu u_2 = (\phi_{is_1,\tau_1} \phi_{is_2, \pi_2})\bu u.
\end{align*}
So we have $u_{12}= u_{21}$.
It follows from Theorem \ref{aw} that
$E_j= \phi_{is_j,\pi_j} E\phi_{is_j, \ti\pi_j}^{-1}$ is the frame of $u_j$ with $E_j(0,0,\l)= \I$. Let $\o_1(x,t)$ and $\o_2(x,t)$ be Hermitian projections onto $E_2(x,t,-is_1)(\Im\tau_1)$ and $E_1(x,t,-is_2)(\Im\tau_2)$ respectively. It follows from Theorem \ref{aw} that we have 
\begin{align*}
E_{12}&= \phi_{is_2,\tau_2} E_1\phi_{is_2, \o_2}^{-1}=  \phi_{is_2,\tau_2}\phi_{is_1, \pi_1}E\phi^{-1}_{is_1,\ti\pi_1}\phi_{is_2, \o_2}^{-1},\\
E_{21}&= \phi_{is_1,\tau_1} E_2\phi_{is_1, \o_1}^{-1}=  \phi_{is_1,\tau_1}\phi_{is_2, \pi_2}E\phi^{-1}_{is_2,\ti\pi_2}\phi_{is_1, \o_1}^{-1}
\end{align*}
are the frame for $u_{12}$ and $u_{21}$ with $E_{12}(0,0,\l)= E_{21}(0,0,\l)=\I$. Since $u_{12}= u_{21}$, we have $E_{12}= E_{21}$. This implies that 
$$\phi_{is_2,\o_2}\phi_{is_1, \ti\pi_1}= \phi_{is_1, \o_1} \phi_{is_2, \ti\pi_2}.$$
 Proposition \ref{av1}(3) implies that $\o_i= \ti\tau_i$ for $i=1,2$.  
\end{proof}

\beg ({\bf Soliton solutions of \eqref{al}})\label{ci}\ 

\ni Note that $u=0$ is a solution of the $(2j-1)$-th vmKdV$_n$ flow \eqref{al} and $E(x,t,\l)=  \exp(a\l x+ a\l^{2j-1} t)$. So 
$$E(x,t,\l)= e^{a\l x+ a\l^3 t}=\bpm  \cos(\l x+ \l^3 t) &-\sin(\l x +\l^3 t) &0\\  sin(\l x+ \l^3 t) & \cos(\l x+ \l^3 t)& 0\\ 0&0 & \I_{n-1}\epm$$
is the frame of $u=0$ with $E(0,0,\l)= \I_{n+1}$. Let $s_1, \ldots, s_k\in \R\bh 0$ satisfying $s_i^2\not= s_j^2$ for all $1\leq i\not= j\leq k$, and $\pi_1, \ldots, \pi_k$ Hermitian projections satisfying \eqref{au}.
 We apply B\"acklund Transform Theorem \ref{aw} to $u=0$ with $E$ and $\phi_{is_j, \pi_j}$ to obtain explicit 1-soliton soluton $u_j= \phi_{is_j, \pi_j}\bu 0$ of \eqref{al}  and its frame
$$E_j(x,t,\l)= \phi_{is_j, \pi_j} E(x,t,\l) \phi^{-1}_{is_j, \ti\pi_j(x,t)} (\l).$$
 We apply the Permutability formula (Theorem \ref{bq}) to $u_i, u_j$ to obtain explicit 2-soliton solutions $u_{ij}$ algebraically from $u_i$ and $u_j$, etc. 
\eeg

\esub

\ss
\bsub {\bf A Poisson structure for the vmKdV$_n$ hierarchy}\

We use a known general method (cf. \cite{DS84}, \cite{Ter97}) to construct a Poisson structure for the vmKdV$_n$ flow and its commuting Hamiltonians. 
 
Identify $\R^{n-1}$ as $Y$ via \eqref{bs}. 
Then 
$$\li\Psi(z_1), \Psi(z_2)\ri:= -\frac{1}{2}\oint \tr(\Psi(z_1)\Psi(z_2))\rd x= \oint z_1^t z_2\rd x =\li z_1, z_2\ri.$$
So $\li\, ,\ri$ is the standard $L^2$ inner product  on $C^\infty(S^1, \R^{n-1})$.

Given $F:C^\infty(S^1,\R^n)\to \R$ and $k\in C^\infty(S^1,\R^{n-1})$, $\K F(k)$ is defined by 
$$\rd F_k(z)= \li \K F(k), z\ri$$
for all $z\in C^\infty(S^1, \R^{n-1})$. 

 Given $k, z\in C^\infty(S^1, \R^{n-1})$, let 
\beq\label{ck}
\Xi_k(z)= z_x-\xi k, \quad {\rm where\,\,}  \xi_x= k z^t - zk^t.
\eeq
Then 
\beq\label{ck1}
\{ H_1, H_2\}(k)=\li \Xi_k(\K H_1(k)), \K H_2(k)\ri
\eeq
is a Poisson structure on $C^\infty(S^1, \R^{n-1})$ and the Hamiltonian flow for a functional $H$ is 
\beq\label{ck2}
k_t= \Xi_k(\K H(k))].
\eeq

\bthm\label{bh} (\cite{Ter97}) Let $a$, $u=\Psi(k)$, and $Q_i(u)$ be as given in Theorem \ref{ak}, and 
$y_{2j-1}(k)$ the $(21)$ block of $Q_{2j-1}(u)$ as defined by \eqref{ake}. Let
\beq\label{bh1} 
F_{2j-1}(k)= -\frac{1}{2j-1} \oint y_{2j-1}(k)\rd x.
\eeq
  Then we have
\ben
\item $\K F_{2j-1}(k)= z_{2j-2}(k)$,
\item the Hamiltonian flow for $F_{2j-1}$ is the $(2j-1)$-th vmKdV$_n$ flow \eqref{al1},
\item $\{F_{2j-1}, F_{2\ell-1}\}=0$,
\een
where $z_{2j-2}(k)$ is the $(32)$-block of $Q_{2j-2}(u)$ as defined by \eqref{ako}.
\ethm

For example, when $j=2$, we have 
$$F_3(k)= -\frac{1}{3} \oint y_3(k)\rd x= -\frac{1}{3} \oint k^t k_{xx} -\frac{||k_x||^2}{2} + \frac{3}{8} ||k||^4 \rd x.$$
 $\K F_3(u)=z_2(k) = -(k_{xx}+ \frac{1}{2} ||k||^2 k)$. By \eqref{ax4}, we have $(\xi_2)_x= k z_2^t - z_2 k^t$. This implies that $\Xi_k(\K F_3(k))= (z_2(k))_x - \xi_2 k$, which is the vmKdV$_n$, the third flow. 
\esub

It follows from  $\{F_{2j-1}, F_{2\ell-1}\}=0$ that we have 

\bprop\label{bx}\

\ben
\item[(i)] $F_{2j-1}$ is a constant of the motion for the $(2\ell-1)$-th vmKdV$_n$ flow, 
\item[(ii)] the Hamiltonian flows for $F_{2j-1}$ and $F_{2\ell-1}$ commute, i.e., the $(2j-1)$-th and $(2\ell-1)$-th vmKdV$_n$ flows commute. 
\een
\eprop 

\bs
\section{The $(2j-1)$-th Airy Curve flow}\label{ay}

In this section, we
\ben
\item construct a parallel frame for a solution of the order $(2j-1)$ curve  flow in $\R^n$ for all $j\geq 1$ so that the principal curvatures satisfy the  $(2j-1)$-th vmKdV$_n$-flow  \eqref{al},
\item give a Poisson structure for the space $\cm_n$ so that the $(2j-1)$-th Airy curve flow is Hamiltonian, 
\item these curve Airy curve flows commute, and admit a sequence of commuting Hamiltonians. 
\een

 \bprop\label{ao} Let $g:\R\to \R^n$ be parameterized by arc-length, and $g=(e_1, \ldots, e_n)$ a parallel frame along $\g$,  $k_i$ the  principal curvature of $\g$ with respect to $e_{i+1}$ for $1\leq i\leq n-1$, and $k=(k_1, \ldots, k_{n-1})^t$. Let $u=\Psi(k)$ as defined by \eqref{bs}, and  $y_{2j-3}(k), \eta_{2j-3}(k)$  the $(21)$ and $(31)$ block of $Q_{2j-3}(u)$ defined by \eqref{ake}. Then
 \beq\label{ao1}
  \g_t= y_{2j-3}e_1+ (e_2,\ldots, e_n)\eta_{2j-3}
  \eeq
  is a well-defined arc-length preserving curve flow on $\R^n$. 
  \eprop 
  
  \begin{proof} First we claim that the right hand side of \eqref{ao1} does not depend on the choice of parallel frames.
  If $g_1$ is another parallel frame along $\g$, then there exists a constant $C\in SO(n-1)$ such that $g_1= g\diag(1,C)$ and the corresponding curvature 
  \beq\label{br1}
  \ti k= C^{-1} k.
  \eeq 
  So $\ti u=\Psi(\ti k)= h^{-1}u h$, where $h=\diag(1,1,C)$. Write 
  $\ti g= (\ti e_1, \ldots, \ti e_n)$. Then we have
  \beq\label{br2}
  \ti e_1= e_1, \quad (\ti e_2, \ldots, \ti e_n)= (e_2, \ldots, e_n) C
  \eeq
   Since $h$ is a constant, it follows from the uniqueness of Theorem \ref{ak} that we have 
  $$Q(\ti u, \l)= Q(h^{-1} uh)= hQ(u,\l) h^{-1}.$$ This implies that $Q_{2j-3}(\ti u)= h^{-1}Q_{2j-2}(u)h$. Recall that $h=\diag(1,1,C)$. So we have 
  \begin{align}
   & y_{2j-3}(\ti u)= y_{2j-3}(u), \label{br3} \\
  & \eta_{2j-3}(\ti k) = C^{-1} \eta_{2j-3}(k) \label{br4}.
  \end{align}
  It follows from \eqref{br1}, \eqref{br2}, and \eqref{br3} that $\ti y_{2j-3}(\ti u)\ti e_1+ (\ti e_2, \ldots, \ti e_n)\eta_{2j-3}(\ti u)$ is equal to $y_{2j-3}( u) e_1+ ( e_2, \ldots,  e_n)\eta_{2j-3}(u)$. This proves the claim. Hence \eqref{ao1} defines a curve flow on $\R^n$.  
  
  By \eqref{ax2}, we have $(y_{2j-3})_x= k^t \eta_{2j-3}$. It now follows from Lemma \ref{ada} that \eqref{ao1} preserves the arc-length parameter. 
  \end{proof} 
  
\bdefn
  We call \eqref{ao1} the {\it $(2j-1)$-th Airy curve flow on $\R^n$\/}. 
  \edefn
  
  \beg
  Since $$(\K^\perp_{e_1})^\ell H= \sum_{i=1}^{n-1} (k_i)_x^{(\ell)} e_{i+1},$$
we can write \eqref{ao1} in terms of $(\K_{e_1}^\perp)^i H$ with $i\geq 0$. For example,  when $j=1$, we have $Q_1(u)= a$. So \eqref{ao1} is the translation flow,
 $$\g_t= e_1= \g_x.$$ 
 For $j=2$, we have 
 $$Q_3(u)= \bpm 0&\frac{||k||^2}{2} & k^t_x\\ -\frac{||k||^2}{2} &0&0\\ -k_x&0&0 \epm,$$
 and the third Airy curve flow \eqref{ao1} is geometric Airy flow \eqref{ba1}. 

The fifth Airy curve flow is $\g_t= y_5(k) e_1+ (e_2, \ldots, e_n) \eta_5(k)$. Use   \eqref{ak6} to rewrite the fifth in terms of $H, \K^\perp_{e_1}H, \cdots, (\K^\perp_{e_1})^{(3)}H$ as follows:
\begin{align*}
\g_t &= (\li H, (\K_{e_1}^\perp)^2 H\ri -\frac{1}{2} ||\K_{e_1}^\perp||^2 + \frac{3}{2} ||H||^4)e_1\\
&\quad + ((\K_{e_1}^\perp)^3 H + \frac{3}{8} ||H||^2 \K_{e_1}^\perp H)
\end{align*}
\eeg
  
 \bthm \label{ao2}If $\g:\R^2\to \R^n$ is a solution of \eqref{ao1}, then there exists $g=(e_1, \ldots, e_n):\R^2\to SO(n)$  such that $g(\cdot, t)$ is a parallel frame along $\g(\cdot, t)$ for each $t$, and $g$ satisfies
 \beq\label{ao3}
 \bca g^{-1}g_x=\bpm 0&-k^t\\ k&0\epm,\\  g^{-1}g_t= \bpm 0& -z_{2j-2}^t(k)\\ z_{2j-2}(k) & \xi_{2j-2}(k)\epm, \eca
 \eeq
 where $z_{2j-2}$ and $\xi_{2j-2}$ are the $(32)$ and $(33)$ blocks of $Q_{2j-1}(u)$ defined by \eqref{ake}. Then
 \ben
 \item  $u=\Psi(k)$ (defined by \eqref{bs}) is a solution of \eqref{al},
  \item if $g_1$ is another map satisfies \eqref{ao3} with $\ti k$ then there exists $c\in SO(n-1)$ such that $g_1= g\diag(1, c)$ and $\ti k= c^{-1}k$. 
  \een
 \ethm
 
 \begin{proof} The proof of this Theorem is similar to that of Theorem \ref{af}. We will use same notations. (2) can be proved exactly the same way. For (1),  
let $h=(e_1, v_2, \ldots, v_n):\R^2\to SO(n)$ such that $h(\cdot, t)$ is a parallel frame along $\g(\cdot, t)$, and $\mu_1, \ldots, \mu_{n-1}$ the corresponding principal curvatures. Then 
$$A:= h^{-1}h_x= \bpm 0& -\mu^t\\ \mu &0\epm.$$
Let $$\hat u=\Psi(\mu),$$ where $\Psi$ is defined by \eqref{bs}. 
Let $B=(b_{ij}):= h^{-1} h_t$. Since $\g$ satisfies \eqref{ao1}, we compute directly to see that 
$$b_{i+1, 1}= \li (e_1)_t, v_{i+1}\ri = \li (\g_t)_x, v_{i+1}\ri = y_{2j-3}(\mu) \mu_i + ((\eta_{2j-3}(\mu))_i)_x.$$
By \eqref{ax1},  $b_{i+1, 1}= (z_{2j}(\mu))_i$, where $z_{2j}(\mu)$ is the $(32)$ block of $Q_{2j}(\hat u)$ as in \eqref{ako}.  So we can write 
$B=\bpm 0 &-z_{2j}^t\\ z_{2j} & S\epm$, for some $S=S(k)\in so(n-1)$. Compare the $(22)$ block of $A_t= B+[A, B]$ to see that $S_x= \mu z_{2j}^t(\mu)- z_{2j}(\mu)\mu^t$.  It follows from \eqref{ax4}, $S_x= (\xi_{2j}(\mu))_x$. 
So there exists $\eta(t)\in so(n-1)$ such that  $S(x,t)= \xi_{2j}(\mu(x,t)) + \eta(t)$. Let $f(t)\in SO(n-1)$ satisfying $f_t= -\eta(t)f$. Then $g=h\diag(1, f)$ is a parallel frame, $k = f^{-1}\mu$ is its principal curvature, and $g$ satisfies \eqref{ao3}.

It follows from Lemma \ref{bb} that $k$ is a solution of \eqref{al1}. 
\end{proof}

Conversely, given a solution of the $(2j-1)$-th vmKdV$_n$ flow \eqref{al1} we can construct a solution of the $(2j-1)$-th Airy curve flow \eqref{ao1}. 

\bprop Let $k$ be a solution of $(2j-1)$-th vmKdV$_n$ flow \eqref{al1}. Then:
\ben
\item There is $g=(e_1, \ldots, e_n):\R^2\to SO(n)$ satisfying \eqref{ao3}. 
\item 
$$\g(x,t)= \int_0^t Z_{2j-3}(0, t_1)\rd t_1 +\int_0^x e_1(s,t)\rd s$$
is a solution of \eqref{ao1}, where 
$$Z_{2j-3}= y_{2j-3}e_1 +(e_2, \ldots, e_n) \eta_{2j-3},$$  $y_{2j-3}(k),$  $\eta_{2j-3}(k)$ are the $(21)$ and $(31)$ blocks of $Q_{2j-3}(u)$ as in \eqref{ake}, and $u=\Psi(k)$. 
\een
\eprop

\begin{proof}
(1) follows from Proposition \ref{cg}.

(2) The proof is similar to that of Proposition \ref{cc}. Note that $\g_x= e_1$.  The second equation of \eqref{ao3} gives
$$(e_1)_t= (e_2, \ldots, e_n) z_{2j-2}.$$
Use the first equation of \eqref{ao3} (which is \eqref{ab1}) and \eqref{ax2} to compute directly to see that 
$$(e_2, \ldots, e_n) z_{2j-2}=(Z_{2j-3})_x.$$ 
Hence \begin{align*}
\g_t &= Z_{2j-3}(0,t) +\int_0^x (e_2, \ldots, e_n)z_{2j-2} \rd s\\
&= Z_{2j-3}(0,t)+ \int_0^x (Z_{2j-3})_s \rd s= Z_{2j-3},
\end{align*} which is the $(2j-1)$-th Airy curve flow \eqref{ao1}. 
\end{proof}

 Cauchy problem on the line for the $(2j-1)$-th Airy curve flow \eqref{ao1} can be proved in the same way as for the geometric Airy curve flow, and similarly for the periodic case. 
 
 \bsub{\bf Poisson structure on $\cm_n$ and commuting Airy curve flows}\

Let $\cm_n$ be as in \eqref{cj}. We have seen that if both $g, \ti g$ are  parallel frames along $\g\in \cm_n$ and $k, \ti k$ are the corresponding principal curvatures, then there exists a constant $c\in SO(n-1)$ such that $\ti g= gc$ and $\ti k= c^{-1}k$.    Note that the group $SO(n-1)$ acts on $C^\infty(S^1, \R^{n-1})$ by $c\dot k= c^{-1}k$. So 
$$\G:\cm_n\to C^\infty(S^1, \R^{n-1})/SO(n-1)$$
defined by $\G(\g)= $ the $SO(n-1)$-orbit of $k$ is well-defined, where $k$ is a the curvature along some parallel frame $g$ along $\g$.

\blem  Let $k, z\in C^\infty(S^1, \R^{n-1})$, $c\in SO(n-1)$, and $\Xi$ the operator defined by \eqref{ck2}. Then 
\beq\label{cm}
\Xi_{c^{-1}k} (c^{-1} z)= c^{-1} \xi_k (z).
\eeq
\elem

\begin{proof}
By \eqref{ck2}, we have $\Xi_k(z)= z_x-\xi k$, where $\xi_x= kz^t- zk^t$. Let $\ti k= c^{-1}k$, $\ti z= c^{-1} z$, and $\ti \xi= c^{-1}\xi c$. A direct computation implies that $\Xi_{\ti k}= \ti k_x - \ti\xi \ti k$ and $\ti\xi_x= \ti k\ti z^t - \ti z \ti k^t$. This proves \eqref{cm}.  
\end{proof}

\bcor If $H_1, H_2$ are functionals on $C^\infty(S^1, \R^{n-1})$ invariant under the action of $SO(n-1)$, then $\{H_1, H_2\}$ is also invariant under the action of $SO(n-1)$. 
\ecor

Formula \eqref{br3} in the proof of Proposition \ref{ao} implies the following:

\bprop\label{cb} If $H:C^\infty(S^1, \R^{k-1})\to \R$ is a functional, then $\hat H:\cm_n\to \R$ defined by $\hat H= H(\G(\g))$ 
is a well-defined functional.
\eprop

Hence $\{\, ,\}$ defined by \eqref{ck1} can be viewed as a Poisson structure on the orbit space  $C^\infty(S^1, \R^{n-1})/SO(n-1)$.  Let $\{\, ,\}^\wedge$ denote the pull back  of $\{\, ,\}$ on $\cm_n$ by the map $\G$, i.e., 
$$
\{\hat H_1, \hat H_2\}^\wedge(\g)= \{H_1, H_2\}(k),
$$
where $\hat H_i= H_i(\G(\g))$.

\bprop\label{cn}  Let $H$ be a functional on $C^\infty(S^1, \R^{n-1})$, and $\K H(k)= z= (z_1, \ldots, z_{n-1})^t$. Then the Hamiltonian vector field for $\hat H$ with respect to $\{\,, \}^\wedge$ is
$$\d \g= A_0e_1+ \sum_{i=1}^{n-1} A_i e_{i+1}$$
where $A_0, A_1, \ldots, A_{n-1}$ satisfy
\beq\label{cn1}
\bca (A_0)_x= \sum_{i=1}^{n-1} k_i A_i,\\ (A_i)_x+ A_0 k_i= z_i,
\eca
\eeq
$g=(e_1, \ldots, e_n)$ is a parallel frame along $\g$ and $k$ is the principal curvature with respect to $g$.
\eprop

\begin{proof} Let  $\d k= \Xi_k(\K H(k))$ denote the Hamiltonian vector field of $H$ with respect to $\{\, ,\}$.  
Since the Poisson structure $\{\, ,\}^\wedge$ is the pull back of $\{\, ,\}$ by $\G$, we have $\rd \G_\g(\d \g)= \d k$. 

Use $g^{-1}g_x= \bpm 0& -k^t\\ k&0\epm$ and a simple computation to see that 
\beq\label{cn2}
[\p_x+\bpm 0& -k^t\\ k&0\epm, g^{-1}\d g]=\bpm 0& -(\d k)^t\\ \d k&0\epm,
\eeq 
 where $\d g$ the variation of parallel frames when we vary $\g$. Since $g\in SO(n)$, $g^{-1}\d g$ is $so(n)$-valued. 
Write $$g^{-1}\d g= \bpm 0& -\eta^t\\ \eta & \zeta\epm.$$
Then \eqref{cn2} implies that $(\zeta)_x= k z^t - z^tk$ and
$\Xi_k(\eta)=\d k= \Xi_k(\K H(k))$. Hence $\eta= z$. 
So we have 
\beq\label{cn3}
g^{-1} \d g= \bpm 0&-z^t\\ z& \zeta\epm,
\eeq
where $z= \K H(u)$. 
 The first column of \eqref{cn3} implies that $$\d e_1= g \bpm 0\\ z\epm.$$
Compute directly to see that 
\begin{align*}
\d e_1&= (\d \g)_x = (A_0e_1+ \sum_{i=1}^{n-1} A_i e_{i+1})_x\\
&= ((A_0)_x- \sum_{i=1}^{n-1} A_i k_i) e_1 + \sum_{i=1}^{n-1} A_0 k_i+ (A_i)_x) e_{i+1} = \sum_{i=1}^{n-1} z_i e_{i+1}.
\end{align*}
This proves \eqref{cn1}.  
\end{proof}

\bprop Let $F_{2j-1}$ be as defined by \eqref{bh1}. Then the Hamiltonian flow for $\hat F_{2j-1}$ on $\cm_n$ is the $(2j-1)$-th Airy curve flow \eqref{ao1}.
\eprop

\begin{proof} Let $y_{2j-3}, \eta_{2j-3}$ be as in \eqref{ake}, and $z_{2j-2}, \xi_{2j-2}$ as in \eqref{ako}. 
By Theorem \ref{bh}, we have $\K F_{2j-1}(k)= z_{2j-2}$. It follows from  \eqref{ax1} and \eqref{ax2} that $A_0= y_{2j-3}$ and $A_i= \eta_{2j-3}$ satisfy \eqref{cn1}. This proves the Proposition. 
\end{proof}

Since $\{ \hat F_{2j-1}, \hat F_{2\ell-1}\}^\wedge= \{ F_{2j-1}, F_{2\ell-1}\}=0$, we have the following:

\bcor\

\ben
\item The $(2j-1)$-th and the $(2\ell-1)$-th Airy curve flows commute. 
\item If $\g(x,t)$ is a solution of the $(2j-1)$-th Airy flow, then $\hat F_{2\ell-1}(\g(\cdot, t))$ is constant in $t$.  
\een
\ecor

\esub 
   
 \bs
 \section{BTs for the geometric Airy flow on $\R^n$}\label{az2}
 
 In this section, we
\ben
\item construct  solutions of the $(2j-1)$-th Airy curve flow \eqref{ao1} from frames of $(2j-1)$-th vmKdV$_n$ \eqref{al},
\item give BT for the $(2j-1)$-th Airy curve flow \eqref{ao1},
\item write down explicit 1-soliton solutions for the geometric Airy flow on $\R^2$. 
\een 

The following Theorem shows that the construction given by Pohlmeyer (\cite{Poh}) and Sym (\cite{Sym})   in soliton theory gives solutions of the $(2j-1)$-th Airy curve flow \eqref{ao1} from frames of solutions of vmKdV$_n$ flow \eqref{al}.

\bthm\label{be} Let $u=\Psi(k)$ be a solution of the $(2j-1)$-th vmKdV$_n$ flow \eqref{al}, $E(x,t,\l)$ the frame of $u$ with $E(0,0,\l)=\I_{n+1}$. Let
$$\zeta(x,t):= E_\l E^{-1}(x,t,0).$$
Then 
\ben
\item $\zeta= \bpm 0& -\g^t\\ \g&0\epm$ for some $\g:\R^2\to \R^n$,
\item $E(x,t,0)= \diag(1, g(x,t))$ for some $g:\R^2\to SO(n)$, $g(\cdot, t)$ is a parallel frame of $\g(\cdot,t)$ and  $k(x,t)$ is the principal curvature of $\g(\cdot, t)$ with respect to $g(\cdot, t)$. 
\item $\g$ is a solution of the $(2j-1)$-th Airy curve flow \eqref{ao1}. 
\een
\ethm

\begin{proof} 
(1) Note that the fixed point set of $\tau$ is $SO(n+1)$ and the connected component of the fixed point set of $\tau$ on $SO(n+1)$ is $K= \{\diag(1, f)\n   f\in SO(n)\}$. 
The $(\tau, \sigma)$-reality condition and $E(0,0,0)=\I_{n+1}$ implies that $E(x,t,0)\in K$. So  
$$f(x,t):=E(x,t,0)=  \diag(1, g(x,t))$$
for some $g(x,t)\in SO(n)$. The $(\tau,\sigma)$-reality condition also implies that
$$\zeta:= E_\l E^{-1}\n_{\l=0}$$ lies in $\cp$.  So there exist $\g:\R^2\to \R^n$ such that
\begin{align*}
&f(x,t):=E(x,t,0)=  \diag(1, g(z,t)), \\
&\zeta=\bpm 0&-\g^t\\ \g&0\epm.
\end{align*}

(2) Since $E$ is a solution of \eqref{an1}, we have 
\beq\label{be1}
\bca E_x= E(a\l + u), \\ E_t= E(a\l^{2j-1} + u \l^{2j-2} +\cdots + Q_{2j-3}(u)\l + Q_{2j-2}(u)),\eca
\eeq
Use \eqref{be1},  $E(x,t,0) = f= \diag(1, g)$, and a direct computation to get
\begin{align}
\zeta_x&=\bpm 0&-\g^t_x\\ \g_x &0\epm= f a f^{-1}= \diag(1, g) a \diag(1, g^{-1})= \bpm 0& -e_1^t\\ e_1 &0\epm, \label{be2}\\
\zeta_t&=\bpm 0&-\g^t\\ \g_t&0\epm=  (E_t)_\l E^{-1}- E_\l E^{-1} E_t E\big|_{\l=0}= fQ_{2j-3}(u) f^{-1}  \\
& =\bpm 1&0\\ 0& g\epm \bpm 0& -y_{2j-3}& -\eta_{2j-3}^t\\ y_{2j-3} &0&0 \\
 \eta_{2j-3}&0&0  \epm \bpm 1& 0\\ 0& g^{-1}\epm. \label{be3}
\end{align}
\eqref{be2} implies that $\g_x= e_1$ the first column of $g$. The first equation of \eqref{be1} at $\l=0$ gives $g^{-1}g_x= \bpm 0&-k^t\\ k&0\epm$. So $g(\cdot, t)$ is a parallel frame and $k(\cdot, t)$ the principal curvature of $\g(\cdot, t)$. 

(3) Write $g=(e_1, \ldots, e_n)$. (3) follows from \eqref{be3}. 
\end{proof}

\bcor\label{cd} Let $\g$ be a solution of \eqref{ao1}, and $g, k$ as in Theorem \ref{ao2}. Let $E$ be a frame of the solution $u=\Psi(k)$ with $E(x,t,0)= g(0,0)$. Then there is a constant $p_0\in \R^n$ such that $E_\l E^{-1}= \bpm 0 & -\g^t-p_0^t\\ \g+ p_0 &0\epm$. 
\ecor

\begin{proof}
Theorem \ref{ao2} implies that $E_\l E^{-1}= \bpm 0&-\hat \g^t\\ \hat \g &0\epm$, $\hat\g$ is a solution of \eqref{ao1}, $g(x,t)= E(x,t,0)$ is a parallel frame for $\hat \g$. But $g$ is also a parallel frame for $\g$. So $\g_x= \hat \g_x$.  Both $\hat\g$ and $\g$ are solutions of \eqref{ao1} implies that 
$$\hat \g_t=\g_t= y_{2j-3}(k)e_1 + (e_2,\ldots, e_n)\eta_{2j-3}(k).$$
Hence $\hat \g-\g$ is a constant $p_0$. 
\end{proof}

\bs

\bthm[{\bf B\"acklund transformation for Airy curve flows}]\label{bg}\

\ni
Let $\g$ be a solution of the $(2j-1)$-th Airy curve flow \eqref{ao1}, and $g(\cdot, t)$ parallel frame and $k(\cdot, t)$ the corresponding principal curvature of $\g(\cdot, t)$ as in Theorem \ref{ao2}.  Let $E$ the frame of the solution $u= \Psi(k)$ of the $(2j-1)$-th vmKdV$_n$ \eqref{al} with initial data $E(0,0,\l)= \diag(1, g(0,0))$, $s\in \R\bh 0$,  $\bc\in \R^n$ a unit vector, and $v=\bpm 1\\ i\bc\epm$. Let
$\ti v(x,t):= E(x,t,-is)^{-1} v$.
Then we have the following: 
\ben
\item $E(x,t,0)=\diag(1, g(x,t))$.
\item $\ti v= (y_0, iy_1, \ldots, iy_n)^t$ for some $y_0, y_1, \ldots, y_n: \R^2\to \R$, $y_0$ never vanishes, and $\sum_{i=1}^n y_i^2= y_0^2$.
\item Write $g=(e_1, \ldots, e_n)$, and set
\begin{align}
\ti \bc&=\frac{1}{y_0}(y_1, \ldots, y_n)^t, \label{bga}\\
\ti \g &= -(\I_n- 2\bc\bc^t)(\g-\frac{2}{sy_0} \sum_{i=1}^n y_ie_i), \label{bw}\\
\ti g&= (\I_n-2\bc\bc^t) g (\I_n- 2\ti \bc\ti \bc^t),\label{ch}\\ 
\ti k&= k-\frac{2s}{y_0} (y_2, \ldots, y_n)^t. \label{bg2}
\end{align}
Then $\ti \g$ is a solution of \eqref{ao1}, 
$\ti g$  is a parallel frame of $\ti \g$ with corresponding principal curvatures $\ti k$
\een
\ethm

\begin{proof}\

(1) Since $E(x,t,\l)\in SO(n+1)$ and satisfies the reality condition \eqref{ar1}, $E(x,t,0)= \diag(1, h(x,t))$ for some $h:\R^2\to SO(n)$. It follows from the fact that $E$ is a frame of $u=\Psi(k)$ that $h$ satisfies \eqref{ao3}. So both $g$ and $h$ satisfies the same linear system \eqref{ao3} and $g(0,0)= h(0,0)$. This proves that $g=h$.

(2) follows from Remark \ref{az3}.  Note that $||\ti \bc(x,t)||^2=1$ for all $(x,t)\in \R^2$.

(3) Let $\pi$, $\ti\pi(x,t)$ be the Hermitian projection onto $\C v$ and $\C\ti v(x,t)$ respectively. So we have 
$$\pi= \frac{1}{2} \bpm 1& -i\bc^t\\ i\bc & \bc\bc^t\epm, \quad \ti\pi= \frac{1}{2} \bpm 1 & -i\ti\bc^t\\ i\ti\bc& \bc\bc^t\epm,$$
where $\ti \bc$ is defined by \eqref{bga}. Let $\phi=\phi_{is,\pi}$. Then 
$$\phi(0)= \bpm -1&0\\0& \I_n-2\bc\bc^t\epm, \quad \phi_\l \phi^{-1}(0)= \bpm 0& -\frac{2}{s} \bc^t \\ \frac{2}{s} \bc &0\epm.$$

Apply Theorem \ref{aw} to $u$ with frame $E$ and $\phi_{is, \pi}$ to see that 
\beq\label{bg3}
\ti u= u+ 2is[a, \ti\pi-\overline{\ti\pi}]
\eeq is a new solution of \eqref{al}.  
Write $\ti u=\Psi(\ti k)$, and $u=\Psi(k)$. Then \eqref{bg3} gives \eqref{bg2}.

By Theorem \ref{aw}, $\ti E= \phi_{is,\pi} E \phi^{-1}_{is, \ti\pi}$ is a frame of $\ti u$. 
  Apply Theorem \ref{be} to $\ti E$ to see that  
$$\ti E_\l \ti E^{-1}|_{\l=0}=\bpm 0&-\hat\g^t\\ \hat\g &0\epm$$ and $\hat \g$ is again a solution of \eqref{ao1}. 
By Corollary \ref{cd}, there is a constant $p_0\in \R^n$ such that
$$E_\l E^{-1}= \bpm 0& -(\g+p_0)^t\\ \g+ p_0& 0\epm.$$ 
Use $\ti E= \phi E\ti\phi^{-1}$ to compute $\ti E_\l \ti E^{-1}$ directly to obtain 
$$\hat\g= \frac{2\bc}{s} - A(\g+p_0) + \frac{2}{s} Ag\ti \bc,$$
where $A=\I_n-2\bc\bc^t$. 
Note that $\frac{2\bc}{s} -Ap_0$ is a constant and the $(2j-1)$-th Airy flow \eqref{ao1} is invariant under the translation. Hence $\ti \g$ defined by \eqref{bw} is a solution of \eqref{ao1}.

Theorem \ref{be} implies that $\ti E(x,t,0)=\diag(1, \ti g(x,t))$ and $\ti g(x,t)$ is a parallel frame for $\hat \g$. So it is a parallel frame for $\ti\g$. 
\end{proof}

\beg[{\bf Explicit soliton solutions of \eqref{ba1} in $\R^2$}]\label{bi} \ 

Note that $k=0$ is a solution of the mKdV, and 
$$E(x,t,\l)= \exp(a(\l x+ \l^3 t))$$
is a frame of the solution $u=\Psi(0)=0$ of the third vmKdV$_2$ flow. Since $E_\l E^{-1}|_{\l=0}= ax$, it follows from Theorem \ref{be} that 
$$\g(x,t)= \bpm x\\ 0\epm$$
is a stationary solution of the geometric Airy curve flow on $\R^2$. 

Let   $\pi$ be the Hermitian projection onto $\C v$, where 
 $$v=\bpm 1\\ 0 \\ i\epm = \bpm 1\\ i\bc\epm, \quad \bc= \bpm 0\\ 1\epm. $$
 Note that 
 $$E(x,t,-is)^{-1}= \bpm \cosh D& -i\sinh D  &0\\ i\sinh D & \cosh D &0\\ 0&0&1\epm,$$
 where 
$$D= sx-s^3t.$$
We apply Theorem \ref{bg} to $\g$ with $k=0$ and $E$  and $\phi_{is,\pi}$  and use the same notation as in Theorem \ref{bg}. A direct computation implies that 
$$\ti c(x,t)= (\tanh (sx-s^3t), \sech(sx-s^3 t))^t, \quad A= \diag(-1, 1),$$ and
\beq\label{bi1}
\ti\g(x,t)=-\bpm 0  \\ \frac{2}{s}\epm + \bpm x-\frac{2}{s}\tanh (sx-s^3t)\\ \frac{2}{s}\sech(sx-s^3 t)\epm
\eeq
is a solution of the geometric Airy flow on $\R^2$ and 
$$\ti g= \bpm 1-2\sech^2D& 2\sech D\tanh (sx-s^3t)\\ -2\sech D\tanh D& 1-2\sech^2 D \epm$$ is a parallel frame with principal curvature 
$$\ti k(x,t)= -2s \sech(sx-s^3t),$$
where $D= sx-s^3t$. Note that $\ti k$ is a $1$-soliton solution of the mKdV. 

If $\ti \g$ is a solution of the geometric Airy flow on $\R^2$, then so is $\ti\g+p_0$ for some constant $p_0\in \R^2$. Hence 
$$\g_1(x,t)=\bpm x-\frac{2}{s}\tanh (sx-s^3t)\\ \frac{2}{s}\sech(sx-s^3 t)\epm$$ is also a solution and $\ti k$ is its curvature. Note that 
   $$\g_1(x,t_0)= \g(x-s^2t_0, 0) + (s^2t_0, 0).$$ So $\g_1$ is a self-similar solution of the geometric curve flow on $\R^2$.  Note that the profile of $1$-soliton solution $\g_1$ when $s=1$ of the geometric Airy flow is the plane curve $\g_0(x)= (x-2\tanh x, 2 \sech x)^t$ (see graph of the curve in Figure 1), and $\g_1$ moves to the right but keeps its shape. 
\eeg

    \begin{figure}
  \includegraphics[width=\linewidth]{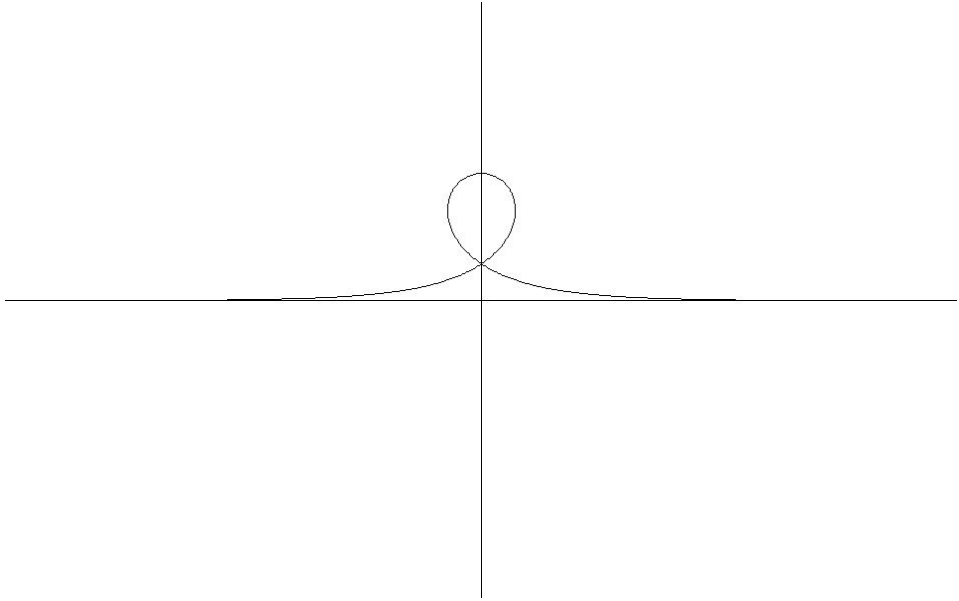}
  \caption{Profile of a 1-soliton solution of  the geometric Airy flow in the plane}
  \label{fig:1-soliton}
\end{figure}

\bs

\beg
Since $k=0$ is a trivial solution of  the mKdV$_n$, 
$E(x,t,\l)= e^{ax+ a^{2j-1}t}$ is the frame of the solution $u=\Psi(0)=0$ of the third flow \eqref{al}. Then $E_\l E^{-1}|_{\l=0}= ax$. By Theorem \ref{be}, $\g(x,t)= (x, 0, \ldots, 0)^t$ is a solution of \eqref{ba1} $\R^n$ with $g(x,t)=\I_n$ as parallel frame and $k=0$ as the corresponding principal curvature. 
Note that
$$E^{-1}(x,t,-is)= \bpm \cosh A& -\sinh A &0\\ \sinh A & \cosh A &0\\ 0&0&\I_{n-2}\epm, \quad A= sx-s^3t.$$
So we can use Theorem \ref{bg} to write down explicit new solutions $\ti \g$ with parallel frame $\ti g$ and principal curvature $\ti k$ as in Example \ref{bi}.  
\eeg

\end{document}